\documentclass[12pt,a4paper]{article}
\usepackage[cp1251]{inputenc}
\usepackage[english]{babel}
\usepackage{graphicx}
\usepackage{amssymb}
\usepackage{amsmath}
\usepackage{amsthm}
\usepackage{comment}
\usepackage{euscript}

\newtheorem{statement}{Statement}
\newtheorem{definition}{Definition}
\newtheorem{theorem}{Theorem}
\newtheorem{lemma}{Lemma}

\newtheorem{prop}{Proposition}
\newenvironment{proposition}{\begin{prop}\rm}{\end{prop}}

\newtheorem{rem}{Remark}

\newtheorem{ex}{Example}

\newtheorem{corollary}{Corollary}

\newcommand{\RP}[0]{\mathbb{RP}}

\newcommand{\be}[0]{\begin{equation}}
\newcommand{\ee}[0]{\end{equation}}
\newcommand{\bez}[0]{\begin{equation*}}
\newcommand{\eez}[0]{\end{equation*}}
\newcommand{\bl}[0]{\begin{lemma}}
\newcommand{\el}[0]{\end{lemma}}

\newcommand{\ep}[0]{$\hspace{\fill} \square$}
\newcommand{\g}[0]{\gamma}

\begin{document}
\noindent
%\textsc{ УДК } 514.774.8

\medskip
\centerline{\Large{Realizability of singular levels of Morse functions}}
\centerline{\Large{by unions of geodesics}}

\medskip
\centerline{\Large{ I.Shnurnikov \footnote{NRU HSE, shnurnikov@yandex.ru}}}

\medskip

{\footnotesize Abstract. We list special graphs of degree 4 with at most 3 vertices (atoms from the theory of integrable hamiltonian systems) which could be represented by a union of closed geodesics on the one of the following surfaces with metric of constant curvature: sphere, projective plane, torus, Klein bottle. 

%Key words:  2--atom, closed geodesics, metric of constant curvature.
}

\bigskip
In the integrable hamiltonian systems theory classifying problems are sometimes solved in terms of neiboughoods of singular levels of functions. A.T.~Fomenko introduced notions of "3--atom" $ $ and "2--atom" $ $ to study the topology of Liuville foliation see \cite{Fomenko 86}, \cite{Fomenko 88} and \cite{Fomenko 89}. Atoms are used in the dinamic systems theory,
%\cite{Fomenko 91 izv an sssr}, \cite{Fomenko 91 funkan}
 groups of symmetries of atoms were stydied in \cite{Kudravzeva 08}.

\begin{definition}
Let $G$ be a two--dimensional compact closed surface, $K$ is a connected finite graph of degree 4, immersed into $G$. The complement $G\setminus K$ is homeomorphic to the union of discs (cells), which could be colored in two colours so that each edge of the graph in incident to cells of different colours. We define atom to be a pair $(G,K)$ up to homeomorphism of pairs (which maps graph into graph).
Advanced atom is a pair $(G,K)$ with fixed coloring of cells into white and black up to homeomorphism which maps cells into the cells of the same color.
\end{definition}
\begin{definition}
An Atom $(G,K)$ is called geodisic if on $G$ there exist a metric of constant curvature such that $K$ could be represented by union of closed geodesics.
\end{definition}

There is one-to-one correspondence between advanced atoms and cell decompositions of two--dimensional closed surfaces with immersions of finite graphs, see \cite[v. 1, \S 2.7.7]{Bolsinov_99}. Graph vertices could have arbitrary degree and decompositions are up to homeomorphism.
For every atom it is possible to find a metric on its surface so that the graph $K$ would be represented by union of closed geodesics. A.T.~Fomenko in 2011 asked what are geodesic atoms with few vertices. 
In denotations of \cite{Bolsinov_99} atom $C_2$ on the sphere, $\widetilde{B}$ on the projective plane, $C_1$ and $E_1$ on the torus are geodesic, atom $B$ on the sphere is't geodesic.
We will say that a way along closed curve (with selfintersections) on a two--dimensional surface changes local orientation if a pair of orthogonal vectors, first of which is tangent to the curve, after continuos deformation along the curve will change the direction of the second vector.\

\begin{proposition}
\label{prop for geodesicity}
If a geodesic atom is represented by the union of closed geodesics $\g_1, \dots, \g_n$, then their intersection points are simple.
A way along closed curve $\g_i$ changes orientation iff the number of intersection points of $\g_i$ with the other geodesics is odd (points of selfintersection of $\g_i$ are not counted ).
\end{proposition}

\begin{proposition}
(a) A set of $n\geq 2$ big circles in the sphere $S^2$ form a geodesic atom iff the circles are in general position.

(b) A set of lines in the projective plane $\RP^2$ form a geodesic atom iff the number of lines is even and they are in general position.
\end{proposition}

\proof
From proposition  \ref{prop for geodesicity} we see that curves are in general position. A way along the line in the projective plane changes the local orientation so by proposition \ref{prop for geodesicity} the  number of lines is even. 
Let us prove by induction on $n$ that these conditions are sufficient. We will form an arrangement by adding circles one after another and change color in all cells on the one side of new circle. For projective plane we fix first line and then change color in all cells on one side of a new line in the disc which is a complement in the plane to the first line of the arrangement.
\ep

\begin{definition}
An admitable pair $(M, \Gamma)$ is a connected two--dimensional compact surface $M$ without boundary and the union $\Gamma$ of a finite set of closed smooth curves on $M$, which intersect (and selfintersect) transversally so that every intersection points belongs to two curves with multiplicities of selfintersection. We alsoo require that $\Gamma$ is connected and contains at least one intersection point.
\end{definition}

Let us consider a closed piece--linear curve $g\subset \Gamma$ with vertices in points of intersections of curves and with edges on the curves on an admitable pair $(M, \Gamma)$. Edges could pass through other intersection points, and in each its own vertex the curve $g$ turns left or right. For every closed curve $g$ let us denote by $u(g)$ the number of points of intersection and selfintersection, which are different from the vertices of $g$ (i.e. $u(g)$ equals to the number of intersection points, in which $g$ doesn't turn).

\begin{statement}
\label{statement orint to number of intersection}
An admitable pair $(M,\Gamma)$ is an atom iff for every nonselfintersecting curve $g \subset \Gamma$ a way along $G$ changes local orientation then and only then, when $u(g)$ is odd.
\end{statement}

\proof
Let pair $(M,\Gamma)$ is an atom with fixed coloring of cells and let $g \subset \Gamma$ be a closed nonselfintersecting curve.
Let us choose a point $P \in g$ and local orientation near point $P$. Let us start a way along $g$ from $P$ and check the color of cell incident to the current edge of way on the right. After passing across the intersection point $u$, which is not a vertex of $g$, the color changes.
% поскольку проход --- сквозной (т.е. переход на противоположное ребро из четырех ребер, исходящих из $u$). 
After passing near intersection point $u$, which is a vertex of $g$, the color doesn't change.  	
%поскольку проход --- поворот (т.е. переход на соседние ребра из четырех ребер, выходящих из $v$). 
%Если число $u(g)$ нечетно(четно), то при обходе вдоль $g$ цвет клетки справа меняется (не меняется) на противоположный, т.е. обход вдоль $g$ меняет (не меняет) локальную ориентацию.

Now we need to color cells in right order to prove the sufficient condition. We will color the regular neighbourhood $\Gamma$ in $M$ by two operations:

(1) If the neighbourhood of an edge $(u,v)$ of graph $\Gamma$ is colored on both sides
%(с одной стороны в белый цвет, с другой стороны --- в черный) 
then we color the neighbourhood  of vertices $v$ (or $u$).

(2) If the neighbourhood of a vertex $v \in \Gamma$ is colored,  then we color correspondingly the neighbourhood  of edges  incident to $v$.

It is easy to see that we will have no contradictions due to sufficiency condition.
%Предположим, что окрашиваемая в некий момент окрестность противоречит с ранее окрашенной, а до этого момента противоречий не было. Тогда противоречащие окрестности соединяются цепочкой правильно окрашенных окрестностей, и все окрестности вместе образуют окрестность замкнутой ломаной $g\subset \Gamma$. Отбрасывая при необходимости в цепочке петли, ломаную $g$ можно сделать без самопересечений. Так мы нашли замкнутую ломаную $g$ без самопересечений, окрестность которой нельзя раскрасить правильным образом. Рассматривая цвета клеток справа при обходе вдоль $g$, аналогично началу доказательства приходим к тому, что для ломаной $g$ не выполняется условие утверждения.
\ep

\begin{corollary}
If a finite set of $n \geq 2$ closed geodesics in the two--dinensional torus $T^2$ is a geodesic atom, then every geodesic contain nonzero even number of intersection points with other geodesics.
\end{corollary}

\begin{proposition}
\label{statement geod in KL^2 orint}
(a) Suppose that $I\subsetneq \g$ is a segment on selfintersecting geodesic $\g$ on the Klein bottle $KL^2$ and the endpoints of $I$ are in the point of selfintersection. Then the way along $I$ changes local orientation.

 \noindent
(b) The boundary of the cell of geodesic atom on the Klein bottle cannot consist of one loop.

 \noindent
(c) Let the point of selfintersection of geodesic  $\g$ divides it on two closed curves $g_1$ and $g_2$ 
%(геодезические с изломом в своих совпадающих концах). 
Then $g_1$ contains the odd number of intersection points with other geodesics and with $g_2$.
\end{proposition}
%\proof
%(а) Рассмотрим локально изометричное универсальное накрытие $p: \R^2 \to KL^2$
% бутылки Клейна. При этом полный прообраз любой точки $X\in KL^2$ состоит из двух подмножеств, для всех точек одного подмножества дифференциал $\mathrm{d}p$ совпадает (касательные пространства отождествлены с $\R^2$), а для точек $x_1$ и $x_2$ разных подмножеств дифференциалы проекций отличаются на отражение $s$ в некоторой прямой $dp_{x_1}= dp_{x_2}\circ s$.
%Поднимем отрезок $I$ до отрезка $AB\subset \R^2$. Если обход вдоль $I$ не меняет локальную ориентацию поверхности, то в точках $A$ и $B$ дифференциалы проекции совпадают. Поэтому проекция отрезка $AB$ на бутылку Клейна даст замкнутую геодезическую (с совпадающими касательными в совпадающих концах), что противоречит выбору отрезка $I$.

%(б) Предположим противное, тогда граница клетки это отрезок замкнутой геодезической без внутренних точек пересечений и самопересечений. По пункту (а) обход вдоль него меняет локальную ориентацию поверхности. Однако этот отрезок (с совпадающими концами) является границей ориентируемой клетки, и в окрестности отрезка можно задать ориентацию клетки --- противоречие.

%(в) Из пункта (а) следует, что обход вдоль отрезка $g_1 \subsetneq \g$ меняет локальную ориентацию. По утверждению \ref{statement orint to number of intersection} (для отрезка $g_1$ в качестве ломаной $g$) получается, что на $g_1$ лежит нечетное число точек пересечения с остальными геодезическими и отрезком $g_2$.
%\ep

\begin{corollary}
\label{corollary KL^2 odd}
One closed selfintersecting geodesic on $KL^2$ with odd number of points of selfintersection cannot form a geodesic atom.
\end{corollary}
%\proof
%Предположим противное. Рассмотрим соответствующую {\it хордовую диаграмму}, т.е. единичную окружность, на которой точки с угловыми координатами $t_1$ и $t_2$ соединены хордой тогда и только тогда, когда точки геодезической с параметрами $t_1$ и $t_2$ совпадают. При этом геодезическая параметризована переменной $t\in [0,2\pi)$. По условию число хорд нечетно. По предложению \ref{statement geod in KL^2 orint}(в) каждая хорда пересекает нечетное число других хорд. Однако это невозможно (число пар пересекающихся хорд в этом случае получается нецелым).  \ep

\begin{theorem}
\label{theorem geodesic atoms}
Among 40 atoms on sphere $S^2$, torus $T^2$, projective plane $\RP^2$ and Klein bottle $KL^2$ which graphs has at most 3 vertices, geodesic are the following:
$C_2$ on $S^2$, $\widetilde{B}$ on $\RP^2$, $C_1$ and $E_1$ on $T^2$, $\widetilde{C}_2,\ \widetilde{D}_2,\ \widetilde{E}_4$ and $\widetilde{G}_3$ on $KL^2$.

\end{theorem}

\proof
Let us give all examples of geodesic atoms. Atom $\widetilde{B}$ on $RP^2$ is formed by 2 projective lines. Atom $C_2$ on $S^2$ is formed by 2 big circles.
%На рис. \ref{pic atom C_1} изображены наборы геодезических на $T^2$, образующие атомы $C_1$ и $E_1$.
Graph of atom $C_1$ consists of a parallel of torus and a geodesic which pass twice along the meridian and once along the parallel. Graph of $E_1$ consists of a parallel, meridian and a diagonal (1,1) of torus.

\begin{figure}[h]
%\vspace*{1cm}
\begin{picture}(450,90)
%бутылка 1
\put(0,0){\circle*{3}}
\put(30,90){\circle*{3}}
\put(120,0){\circle*{3}}
\put(150,90){\circle*{3}}
\put(0,0){\vector(1,0){119}}
\put(150,90){\vector(-1,0){119}}
\put(0,0){\vector(1,3){30}}
\put(120,0){\vector(1,3){30}}
%геодезическая 1
\thicklines
\put(0,0){\line(4,3){120}}
\put(30,0){\line(-4,3){24}}
\put(126,18){\line(-4,3){96}}
%пересечения
\put(15,11.25){\circle*{5}}
\put(75,56.25){\circle*{5}}
%бутылка 2
\put(200,0){\begin{picture}(150,90)
\thinlines
\put(0,0){\circle*{3}}
\put(30,90){\circle*{3}}
\put(120,0){\circle*{3}}
\put(150,90){\circle*{3}}
\put(0,0){\vector(1,0){119}}
\put(150,90){\vector(-1,0){119}}
\put(0,0){\vector(1,3){30}}
\put(120,0){\vector(1,3){30}}
\thicklines
%геодезическая 1
%\put(15,45){\circle*{3}}
\put(15,45){\line(1,0){120}}
%\put(135,45){\circle*{3}}
%геодезическая 2
%\put(75,0){\circle*{3}}
%\put(75,90){\circle*{3}}
\put(75,0){\line(0,1){90}}
%геодезическая 3
\put(15,45){\line(0,-1){45}}
\put(135,45){\line(0,1){45}}
%пересечения
\put(75,45){\circle*{5}}
\put(15,45){\circle*{5}}
\put(135,45){\circle*{5}}
\end{picture}}
\end{picture}
\caption{Atoms $\widetilde{C}_2$ (left) and $\widetilde{D}_2$ on the Klein bottle}
%\vspace*{1cm}
\label{pic atom KL_2 C, D}
\end{figure}
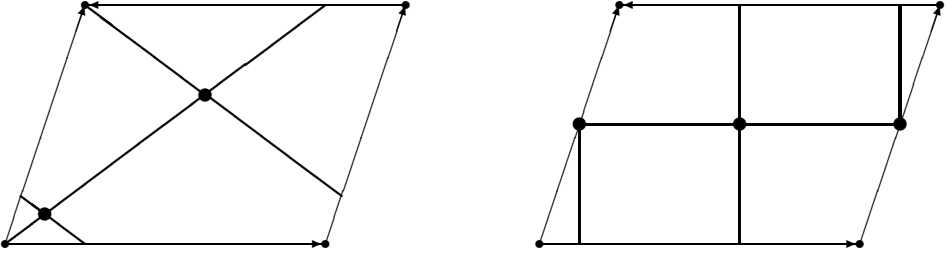

%На рис. 1 и 2 
%\ref{pic atom KL_2 C, D} и \ref{pic atom KL_2 E, G} 
%наборами геодезических на $KL^2$ реализованы атомы $\widetilde{C}_2$, $\widetilde{D}_2$, $\widetilde{E}_4$ и $\widetilde{G}_3$. 

Graph of $\widetilde{C}_2$ consists of one closed geodesics with two points of selfintersecion, graph of atom $\widetilde{D}_2$ consists of a parallel (not changing the orientation) and two meridians (changing the orientation), graph of $\widetilde{E}_4$ consists of parallel and selfintersecting in one point geodesic, graph of $\widetilde{G}_3$ consists of two meridians and selfintersecting in one point geodesic.

\begin{figure}[h]
%\vspace*{1cm}
\begin{picture}(450,90)
%бутылка 1
\put(0,0){\circle*{3}}
\put(30,90){\circle*{3}}
\put(120,0){\circle*{3}}
\put(150,90){\circle*{3}}
\put(0,0){\vector(1,0){119}}
\put(150,90){\vector(-1,0){119}}
\put(0,0){\vector(1,3){30}}
\put(120,0){\vector(1,3){30}}
%геодезическая 1
\thicklines
\put(0,0){\line(2,3){60}}
\put(60,90){\circle*{3}}
\put(90,0){\circle*{3}}
\put(90,0){\line(-2,3){60}}
% геодезическая 2
\put(15,45){\circle*{3}}
\put(135,45){\circle*{3}}
\put(15,45){\line(1,0){120}}
%пересечения
\put(45,67.5){\circle*{5}}
\put(30,45){\circle*{5}}
\put(60,45){\circle*{5}}
%бутылка 2
\put(200,0){\begin{picture}(150,90)
\thinlines
\put(0,0){\circle*{3}}
\put(30,90){\circle*{3}}
\put(120,0){\circle*{3}}
\put(150,90){\circle*{3}}
\put(0,0){\vector(1,0){119}}
\put(150,90){\vector(-1,0){119}}
\put(0,0){\vector(1,3){30}}
\put(120,0){\vector(1,3){30}}
\thicklines
%геодезическая 1
\put(0,0){\line(2,3){60}}
\put(60,90){\circle*{3}}
\put(90,0){\circle*{3}}
\put(90,0){\line(-2,3){60}}
%геодезическая 2
\put(75,0){\circle*{3}}
\put(75,90){\circle*{3}}
\put(75,0){\line(0,1){90}}
%геодезическая 3
\put(15,45){\circle*{3}}
\put(15,45){\line(0,-1){45}}
\put(135,45){\circle*{3}}
\put(135,45){\line(0,1){45}}
%пересечения
\put(45,67.5){\circle*{5}}
\put(75,22.5){\circle*{5}}
\put(15,22.5){\circle*{5}}
\end{picture}}
\end{picture}
\caption{Atoms $\widetilde{E}_4$ (left) and $\widetilde{G}_3$ on the Klein bottle}
%\vspace*{1cm}
\label{pic atom KL_2 E, G}
\end{figure}
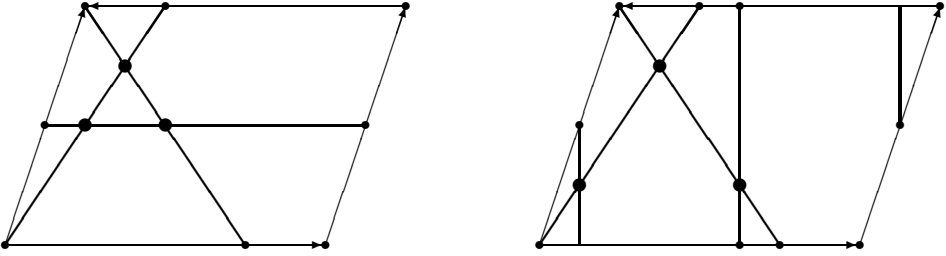

In \cite[v.1, ch.2]{Bolsinov_99} there is a list of atoms with at most three vertices. 
% Для всех из них, кроме реализованных выше и кроме шести атомов на сфере с тремя листами Мебиуса, не попадающих под условие теоремы, мы укажем, почему они не геодезические. 
%Атом $A$ не геодезический, т.к. точка не считается геодезической. 
Atoms
$$
B, D_1, D_2, E_3, F_2, G_1, G_2, G_3, H_1, H_2 \quad \text{on $S^2$ and atoms} \ \ 
\widetilde{C}_1, \widetilde{D}_1, \widetilde{E}_6, \widetilde{F}_5, \widetilde{F}_6, \widetilde{G}_4, \widetilde{G}_5, \widetilde{G}_6, \widetilde{G}_7,\widetilde{H}_3, \widetilde{H}_4 \ \  \text{on}\ \RP^2
$$ atoms on the torus $E_2, F_1$ --- are not geodesic because  closed geodesics on these surfaces with metrics of costant curvature cannot have points of selfintersection. 
The rest atoms are not geodesics because:
$\widetilde{E}_3$ by corollary \ref{corollary KL^2 odd};  $\widetilde{E}_5$, as two simple geodesics on $KL^2$ cannot intersect in three points; $\widetilde{E}_7, \widetilde{F}_7$, as two lines on $\RP^2$ intersects in one point;  $\widetilde{F}_3, \widetilde{G}_2, \widetilde{H}_2$ by propostion \ref{statement geod in KL^2 orint}(b);  $\widetilde{F}_4$by propostion \ref{statement geod in KL^2 orint}(c).
\ep

%\begin{remark}
%Среди всех атомов с не более чем тремя вершинами есть шесть атомов на сфере с тремя листами Мебиуса, о которых еще не известно, будут ли они геодезическими.
%\end{remark}

\noindent
{\bf Acknowledments}. I'm grateful to A.T.~Fomenko for the motivation of this work.

\end{document}